\documentclass[11pt,a4paper,leqno]{article}

\usepackage{color}
\usepackage[x11names]{xcolor}
\usepackage{url}
\usepackage[utf8]{inputenc}
\usepackage[T1]{fontenc}
\usepackage{mathtools,amssymb,mathrsfs}
\usepackage{lmodern}
\usepackage{amsmath}
\usepackage{amsthm}
\numberwithin{equation}{section}
\usepackage{fancybox}
\usepackage{hyperref}

\usepackage{tikz}
\usetikzlibrary{decorations.markings}
\usetikzlibrary{decorations.pathreplacing}

\usepackage{mathtools}

\mathtoolsset{showonlyrefs}

\usepackage{soul}
\def\op{\mathrm{op}}
\def\pdo{\mathrm{pdo}}

\def\v_j{\underline{ {v}^{j} }}
\def\R_j{\mathrm{R}^{j}}

\def\l_j{\widetilde{ \lambda^{ j } }}

\def\o{\omega}

\title{\textbf{On M\'etivier's Lax-Mizohata theorem and extensions to weak defects of hyperbolicity. \\ Part two}}

\date{\today}

\author{Karim Ndoumajoud\thanks{Institut de Math\'ematiques de Jussieu - Paris Rive Gauche UMR CNRS 7586, Sorbonne Universit\'e} \and Benjamin Texier\thanks{Institut Camille Jordan UMR CNRS 5208, Universit\'e Claude Bernard Lyon 1 {\tt texier@math.univ-lyon1.fr}}}

\newtheorem{lem}{Lemma}[section]
\newtheorem{prop}{Proposition}[section]
\newtheorem{theo}{Theorem}

\newtheorem{ass}{Assumption}[section]
\newtheorem{rem}{Remark}[section]

\def\R{\mathbb R}
\def\C{\mathbb C}
\def\S{\mathbb S}
\def\N{\mathbb N} 
\def\op{{\rm op}}
\def\e{\varepsilon}
\def\s{\sigma}
\newcommand{\be}{\begin{equation}}
\newcommand{\ee}{\end{equation}}
\newcommand{\ba}{\begin{aligned}}
\newcommand{\ea}{\end{aligned}}
\renewcommand{\d}{\partial}
\def\a{\alpha}
\def\b{\beta}
\def\g{\gamma}
\def\l{\lambda}

\def\op{{\rm op}}

\begin{document}

\maketitle

\begin{abstract}
We continue our study of initial-value problems for fully nonlinear systems exhibiting strong or weak defects of hyperbolicity. We prove that, regardless of the initial Sobolev regularity, the initial-value problem has no local $H^s$ solutions, for $s > s_0 + d/2,$ if the principal symbol has a strong, or even weak, defect of hyperbolicity, and the purely imaginary eigenvalues of the principal symbol are semi-simple and have constant multiplicity. The index $s_0 > 0$ depends on the severity of the defect of hyperbolicity. These results recover and extend previous work from G. M\'etivier [{\it Remarks on the well posedness of the nonlinear Cauchy problem,} 2005], N.Lerner, Y. Morimoto, C.-J. Xu [{\it Instability of the Cauchy-Kovalevskaya solution for a class of non-linear systems}, 2010] and N. Lerner, T. Nguyen, B. Texier, {\it The onset of instability in first-order systems}, 2018].
\end{abstract}

\section{Introduction}

We continue our study, initiated in \cite{NT1}, of initial-value problems for fully nonlinear systems exhibiting strong or weak defects of hyperbolicity. Our setting is the same as in \cite{NT1}: we consider initial-value problems in $\R^d:$  
\begin{equation} \label{equation reference} 
\d_t u + F(t,x,u,\d_x u)  = 0, \qquad u_{ |t = 0 }  = u_{in} \in H^{ \sigma }( \mathbb{ R }^{ d } )
\end{equation}
where $ t \geq 0$ and $d \geq 1.$ The map $F$ is smooth in its arguments $t \in \R_+,$ $x \in \R^d,$ $u \in \R^N$ and $v \in \R^{Nd},$ and takes values in $\R^N.$ 
The principal symbol is defined as
\be \label{def:A}
A(t,x,u,v,\xi) :=  \xi \cdot \d_4 F(t,x,u,v) =  \sum_{1 \leq j \leq d} \xi_j \d_{v_j} F(t,x,u,v) \quad \in \R^{N \times N}.
\ee

We prove that the initial-value problem has no Sobolev solutions in $H^s(\R^d),$ for $s > s_0 + d/2,$ if the principal symbol $A$ has a strong, or even weak, defect of hyperbolicity, and the purely imaginary eigenvalues of $A$ have constant multiplicity. The integer $s_0 \in \N$ depends on the severity of the defect of hyperbolicity. The results hold regardless of the initial Sobolev regularity. 

Theorem \ref{th:main} below recovers an ellipticity result of G. M\'etivier \cite{Metivier}. Theorem \ref{th:transition} below continues our study of weak defects of hyperbolicity, which builds on work of N. Lerner, Y. Morimoto, and C.-J. Xu, \cite{LMX}, and N. Lerner, T. Nguyen and B. Texier \cite{LNT}. 

Both ill-posedness theorems in this article rely on and extend the ill-posedness results of our previous article \cite{NT1}.   

\section{Assumptions and results}

\subsection{Elliptic systems} 

\begin{ass} \label{ass:ell} For some $(x^0,\xi^0) \in \R^d \times \S^{d-1},$ some $(u^0, v^0) \in \R^N \times \R^{Nd}:$
\begin{itemize}
\item[{\rm (i)}] The spectrum of the principal symbol $A(0,x^0,u^0,v^0,\xi^0)$ defined in \eqref{def:A} is not entirely contained in $\R.$ 
\item[{\rm (ii)}] For all $\xi \in \S^{d-1},$ the real eigenvalues of $A(0,x^0,u^0,v^0,\xi)$ are semi-simple and have constant multiplicity.
\end{itemize}
\end{ass}

The above assumption implies non-existence of Sobolev solutions to the initial-value problem \eqref{equation reference}, in the following sense: 
\begin{theo}\label{th:main}
Under Assumption {\rm \ref{ass:ell},} for any $s >2 + d/2,$ any $\s > 0,$ for some $u_{in} \in H^\sigma(\R^d),$ there is no $T > 0$ and no ball $B_{x^0}$ centered at $x^0$ such that the initial-value problem \eqref{equation reference} has a solution in $C^0([0,T], H^s(B_{x^0})).$  
\end{theo}

Theorem \ref{th:main} recovers an ellipticity theorem of M\'etivier (Theorem 4.5, \cite{Metivier}). Theorem 1 in \cite{NT1} states that under Assumption \ref{ass:ell}(i), the initial-value problem \eqref{equation reference} has no solution if $1 + d/2 < s \leq \s < s + (s - 1 - d/2).$ We see here that the above bound on the initial regularity index $\s$ can be removed under the condition formulated in Assumption \ref{ass:ell}(ii), if in addition the lower bound on $s$ is strengthened to $s > 2 + d/2.$ Note that condition (ii) in Assumption \ref{ass:ell} is global in $\xi \in \S^{d-1}.$ 

 The proof of Theorem \ref{th:main}, given in Section \ref{sec:proof:1} below, shows that under Assumption \ref{ass:bif}, a posited solution $u$ to the initial-value problem \eqref{equation reference} with a priori regularity $H^s$ is in fact $H^{s'}$ over a small time interval, with $s' > s.$ Indeed: by Assumption \ref{ass:ell}, we may decompose the para-linearized system based on \eqref{equation reference} into hyperbolic, positive elliptic, and negative elliptic subsystems. For the hyperbolic subsystem, associated with purely imaginary eigenvalues of the principal symbol evaluated along the posited solution, the initial regularity is propagated. For the positive elliptic equations, associated with eigenvalues of the principal symbol with positive real parts, we observe a forward-in-time regularization effect, and for the negative elliptic equations, we observe symmetrically a backward-in-time regularization effect. Iterating the argument showing a gain of regularity, we arrive at a contradiction with Theorem 1 from \cite{NT1}. This disproves the existence of sufficiently regular Sobolev solutions to the initial-value problem \eqref{equation reference}.  %

 As mentioned in \cite{NT1}, we view the extension of Theorem \ref{th:main} to Gevrey spaces as an important future development. It was our main motivation in writing out a detailed version of M\'etivier's proof. The work of B. Morisse \cite{Mo1,Mo2,Mo3} tackles well-posedness issues in Gevrey spaces.

\subsection{Transition to ellipticity} 

 We extend here Theorem 2 from \cite{NT1}.
We let $P$ be the characteristic polynomial of $A$ evaluated at $u(t,x):$ 
\be \label{def:char}
 P(t,x,\xi,\l) = \mbox{det}\, (A(t,x, u(t,x), \d_x u(t,x),\xi) - \l {\rm Id}\big).
\ee
By the equation \eqref{equation reference}, the Taylor expansion at $t = 0$ of $P$ depends on $F$ and the initial datum $u_{in}.$ Given $U$ an open set in $\R^d \times \S^{d-1},$ we denote ${\mathcal S}_U$ the spectrum at $t = 0$ above $U:$ 
$$ {\mathcal S}_U := \big\{ (x,\xi,\l) \in U \times \C, \quad P(0,x,\xi,\l) = 0 \big\}.$$

The following assumption, based on Hypothesis 1.5 from \cite{LNT}, describes pairs of eigenvalues coalescing and branching out of the real axis at $t = 0,$ in a $C^1$ fashion:

\begin{ass} \label{ass:bif} We assume
\begin{itemize}
\item[{\rm (i)}] {\rm initial hyperbolicity:} for some open set $U,$ we have $\o = (x,\xi,\l) \in {\mathcal S}_U$ only if $\l \in i \R.$
\item[{\rm (ii)}]  For all $\xi \in \S^{d-1},$ the real eigenvalues of $A(0,x^0,u^0,v^0,\xi)$ are semi-simple and have constant multiplicity. 
\item[{\rm (iii)}] {\rm Existence of double eigenvalues:} For some $\o^0 = (x^0,\xi^0,\l^0) \in {\mathcal S}_U,$ we have $\d_\l P(0,\o^0) = 0.$ 
\item[{\rm (iv)}] {\rm $C^1$ bifurcations of the spectrum away from $\R:$} for any $\xi \in \S^{d-1},$ given any $\l \in  \R$ such that $\o = (x^0,\xi,\l) \in {\mathcal S}_U$ and $\d_\l P(0,\o) = 0,$ we have:
\be \label{ineg:bif} \left\{ \begin{aligned}   \d_t P(0,\o') & = 0, \quad \d_\l P(0, \o') = 0, \quad \mbox{for all $\o'$ near $\o$ in ${\mathcal S}_U,$}  \\ 
 \big (\d_{t\l}^2 P(0, \o) \big)^2 & < (\d_t^2 P \d_\l^2 P)(0,\o),  
  \end{aligned}\right.
 \ee    
\end{itemize}
\end{ass}

\begin{theo} \label{th:transition} Under Assumption {\rm \ref{ass:ell},} for any $s > 4 + 3 d/2,$ any $\s > 0,$ for some $u_{in} \in H^\sigma(\R^d),$ there is no $T > 0$ and no ball $B_{x^0}$ centered at $x^0$ such that the initial-value problem \eqref{equation reference} has a solution in $C^0([0,T], H^s(B_{x^0})).$  
\end{theo}

The proof of Theorem \ref{th:transition}, given in Section \ref{sec:weak}, shows a regularization effect and then invokes Theorem 2 from \cite{NT1}. The regularization effect is here degenerate in time, since under Assumption \ref{ass:bif}, the real parts of the eigenvalues of the principal symbol evaluated along the posited solution are $O(t).$ 

In Theorem \ref{th:transition}, the condition on $s$ can be read as $s - 1 - d/2 > 3 + d:$ we need indeed for a posited solution $u$ that $(u, \d_x u)$ belongs to $W^{3 + d,\infty},$ where $3 + d$ is the number of derivatives that comes in the G\r{a}rding inequality (as per, for instance, the proof of Theorem 1.1.26 in \cite{Lerner}). Remark \ref{rem:garding} expounds on this issue.

We refer to \cite{LNT} for examples of systems and data satisfying Assumption \ref{ass:bif}.

\section{Proof of Theorem \ref{th:main}: strong defects of hyperbolicity} \label{sec:proof:1} 

By contradiction, we assume existence of $T > 0,$ a neighborhood $B_{x^0}$ of $x^0$ and $u \in C^0([0,T], H^s(B_{x^0})$ such that $u$ solves the initial-value problem \eqref{equation reference}. We identify $u$ with a map defined in all of $\R^d$ via an extension operator $H^s(B_{x^0}) \to H^s(\R^d).$

If $\s < s + (s - 1 - d/2),$ then Theorem 1 from \cite{NT1} applies, and asserts that there are no solution to the initial-value problem \eqref{equation reference}. Thus we may assume $\s \geq s + (s - 1 - d/2).$ In particular, since $s > 2 + d/2,$ this implies $\s \geq 3 + d/2.$  

The proof goes as follows: first we prove that the posited solution actually belongs to $H^{s + \rho(s)}$ with
\be \label{def:rho}
 \rho(s) =\min (\s - s, \, s - 2 - d/2, \, 1/2) > 0
\ee
 over a smaller time interval $[0, T'],$ with $T' < T.$ 
 
\begin{figure} 
\hspace{1.75cm}
\begin{tikzpicture}


\begin{scope}[shift={(-5,0)}]

\draw[->] (-.3,0) -- (9.5,0) ; 
\draw[->] (0, -.3) -- (0, 4);

\draw (0,4) node[anchor=west] {\footnotesize $\rho$} ;
\draw (9.5,-.2) node[anchor=north] {\footnotesize $s$} ; 


\draw[very thick] (4.64,1) -- (6.07,1) ;

\draw[very thick] (6.07,1) -- (8.5,0) ;

\draw[dotted] (0,3.5) -- (8.5,0) ; 

\draw (8.5,-.2) node[anchor=north] {\footnotesize $\s$} ; 

\filldraw (8.5,0) circle (0.05cm) ; 

\draw (5.2, 3.5) node[anchor=west] {\footnotesize Theorem 1 from \cite{NT1} applies} ; 

\draw (5.2,3) node[anchor=west] {\footnotesize if $s > (\s + 1 + d/2)/2$} ;



\draw[dotted] (0,1) -- (8.5,1) ; 

\draw[dotted] (1,-.5) -- (7,1.97) ;

\draw[very thick] (2.21,0) -- (4.64,1) ; 

\filldraw (2.21,0) circle (0.05cm) ; 

\draw (2.21,0) node[anchor=north] {\footnotesize $2 + d/2$} ; 

\draw[dotted] (6.07,0) -- (6.07,1) ;

\filldraw (0,1) circle (0.05cm) ; 

\draw (0,1) node[anchor=east] {\footnotesize $1/2$} ; 

\draw (0,3.5) node[anchor=east] {\footnotesize $\s$} ;

\filldraw (0,3.5) circle (0.05cm) ; 

\draw (1.4,2.4) node[anchor=south] {\footnotesize $\s - s$} ;  

\draw (7.7,1.5) node[anchor=south] {\footnotesize $s -2 -d/2$} ; 
\draw (6.07, 0 ) node[anchor=north] {\footnotesize $\s - 1/2$} ;

\filldraw (6.07,0) circle(0.05cm) ; 

\draw[red] (5,-.7) -- (5,4) ;

\filldraw (5,0) circle(0.05cm) ;

\draw (5,-.7)  node[anchor=north] {\footnotesize $(\s + 1 + d/2)/2$} ;

\end{scope}

\end{tikzpicture}

\caption{The proof of Theorem \ref{th:main} shows that if a solution to \eqref{equation reference} is in $H^s,$ then on a smaller time interval that solution actually belongs to $H^{s + \rho(s)},$ with $\rho$ given by the thick line on the above picture, in accordance with \eqref{def:rho}. Theorem 1 from \cite{NT1} states that no solution exists in $H^s$ if $s > (\s + 1 + d/2)/2.$ In particular, since $s > 2 + d/2,$ we may assume $\s > \min(3 + d/2, s).$ Given an initial regularity count $s > 2 + d/2,$ observing that $\rho$ is positive and non-decreasing in $[s, (\s + 1 + d/2)/2),$ we see that after applying the regularizing argument a finite number of times we obtain a regularity count that is greater than $(\s + 1 + d/2)/2,$ so that Theorem 1 from \cite{NT1} applies.}  \label{fig1}

\end{figure}
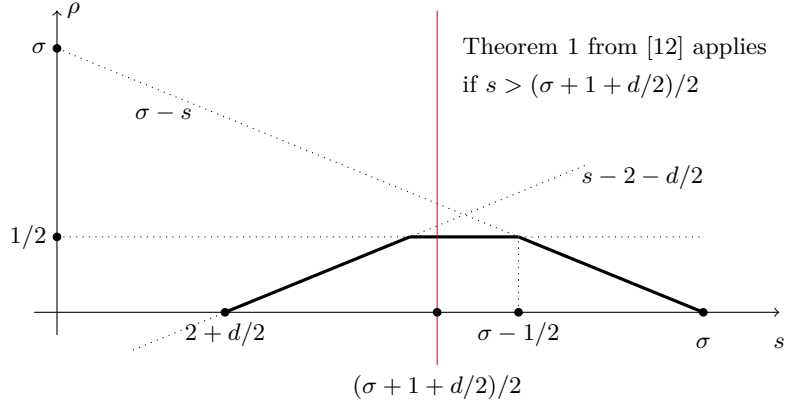

We may use the argument again, so that the solution has regularity $s + \rho(s) + \rho(s + \rho(s)),$ over some $[0, T''],$ with $T'' < T'.$ Since $\rho$ is positive and non-decreasing in $s$ over the relevant $s$ domain (see Figure \ref{fig1}) after a finite number of steps we find that the solution has Sobolev regularity strictly greater than $(\s + 1 + d/2)/2,$ over a small but non-trivial time interval, at which point Theorem 1 from \cite{NT1} brings a contradiction. 

%
%

\subsection{Paralinerization} \label{sec:paralin} 

Our first step is to para-linearize the system. Since $u$ belongs to $H^s$ for $t \in [0, T],$ we have (see for instance Theorem 5.2.4 in \cite{Metivierbis}):
\be \label{para}
 F( t,  x, u, \d_x u) = T_{\d_3 F(u ,  \d_x u)} u + T_{\d_4 F(u , \d_x u)}  \d_x u + R_{para},
\ee
with $R_{para} \in L^\infty([0, T], H^{2 (s -1) - d/2}(\R^d)).$ Above in \eqref{para}, $T_a b$ refers to the para-product of $b$ by $a,$ and $\d_3 F$ refers to the partial derivative of $F$ with respect to its third ($u$) argument. %
We now use notation
$\op(a) = T_a$ to denote the para-differential operator with symbol $a$ (see Appendix \ref{sec:symb}). We denote $\tilde A$ the perturbation of the principal symbol $A$ defined by 
\be \label{def:tildeA}
 \tilde A := \d_3 F + i \xi \cdot \d_4 F.
 \ee 
 Thus \eqref{para} takes the form
$$ F(t,x,u,\d_x u) = \op(\tilde A) u  + R_{para}.$$

\subsection{Localization in a cone in phase space} \label{sec:loc}

Under Assumption \ref{ass:ell}, all frequencies $\xi \in \S^{d-1}$ are such that $A(0,x^0,u^0,v^0,\xi)$ has non-real eigenvalues. Indeed, the set of all $\xi$ such that $A(0,x^0,u^0,v^0,\xi)$ has only real eigenvalues is open by condition (ii) in Assumption \ref{ass:ell}, and closed by continuity of the eigenvalues. Thus it is empty, by condition (i) in Asumption \ref{ass:ell} and connectedness of the sphere.

We single out a frequency $\xi^0 \in \S^{d-1},$ and localize in $(x,\xi)$ in a cone that contains $(x^0,\xi^0) \in \R^d \times \S^{d-1},$ where $x^0 \in \R^d$ is put forward in Assumption \ref{ass:ell}.
In this view, let $U_{x^0}$ be a small neighborhood of $x^0,$ let $V_{\xi^0}$ be a small neighborhood of $\xi^0$ on the sphere,  and 
$$\Omega := \big\{ (x,\xi) \in \R^{2d}, \quad x \in U_{x^0}, \quad \xi/|\xi| \in V_{\xi^0} \quad \mbox{and} \quad |\xi| \geq r \big\} ,$$
for some $r > 0$ which will chosen somewhat large.  
We define $\psi: \R^{2d}_{x,\xi} \to [0,1],$ to be smooth, supported in $\Omega,$ and identically equal to $1$ in a smaller truncated cone $\Omega^\flat,$ defined by 
$$ \Omega^\flat := \{ (x,\xi) \in \R^{2d}, \quad x \in U^\flat_{x^0}, \quad \xi/|\xi| \in V^\flat_{\xi^0} \quad \mbox{and} \quad |\xi| \geq r'\},$$
where $x^0 \in U^\flat_{x^0} \subset U_{x^0},$ and $\xi^0 \in V^\flat_{\xi^0} \subset V_{\xi^0},$ and $r^\flat \geq r.$ We may choose $\psi$ to be a tensor product $\psi(x,\xi) = \psi_1(x) \psi_2(\xi).$ We define similarly $\psi^\flat$ to be another smooth cut-off, supported in $\Omega^\flat$ and identically equal to 1 in yet another, smaller truncated cone. 

We let 
\be \label{def:v} v = \op(\psi) u,\ee so that 
the initial-value problem in $v$ is
\be \label{ivp:v}
 \d_t v + \op(\tilde A) v  = f , \qquad v(0)  = \op(\psi) u_{in},
\ee
where
\be \label{def:f}
 f = \op(\psi) R^{para} - [\op(\psi), \op(\tilde A)] u.
\ee
We observe that the operator in \eqref{ivp:v} is further localized, up to a remainder. Indeed, we have
 $$ \op(\tilde A) v = \op(\psi^\sharp \tilde A) v - f_1,$$
 where
 \be \label{def:f1} f_1 := \op(\tilde A) \big(\op(\psi^\sharp) \op(\psi) - \op(\psi)\big) u + (\op(\psi^\sharp \tilde A) - \op(\tilde A) \op(\psi^\sharp)\big) v.\ee 

Above, the cut-off $\psi^\sharp$ is defined similarly to $\psi,$ on a bigger cone $\Omega^\sharp \supset \Omega,$ and identically equal to $1$ on $\Omega,$ so that $(1 - \psi^\sharp) \psi \equiv 0.$ (The cone $\Omega^\sharp$ is defined just like $\Omega,$ above, in terms of neighborhoods $U^\sharp_{x^0}$ of $x^0$ and $V^\sharp_{\xi^0}$ of $\xi^0,$ and a radius $0 < r^\sharp < r.$ Conditions bearing on $\Omega^\sharp$ will be formulated in Section \ref{sec:spectraldec} and Lemma \ref{lem:gstar}.) 
We let $g = f + f_1,$ 
and the initial-value problem in $v$ now appears as
\be \label{ivp:v:2}
 \d_t v + \op(\psi^\sharp \tilde A ) v  = g,  \qquad 
 v(0)  = \op(\psi) u_{in}.\ee

\subsection{Projections} \label{sec:spectraldec}

 If $T$ and the support of $\psi^\sharp$ are small enough, meaning $r^\sharp$ large enough and $U_{x^0}^\sharp \subset \R^d$ and $V^\sharp_{\xi^0} \subset \S^{d-1}$ small enough, then $\tilde A(t,x,\xi)$ is a small perturbation of $A(t,x,\xi).$ By symmetry (eigenvalues out of $i \R$ occur in pairs) and continuity of the spectrum (see for instance Proposition 1.1 in \cite{Rouche}), the spectrum of $\psi^\sharp \tilde A$ then satisfies Assumption \ref{ass:ell}, for all $t \leq T$ and all $(x,\xi).$ We suppose that these conditions bearing on $T$ and $\psi^\sharp$ are satisfied in the following.   

 We decompose the spectrum into ``hyperbolic'' eigenvalues that belong to $i \R,$ ``negative elliptic'' eigenvalues which have a negative real part, and ``positive elliptic'' eigenvalues, which have a positive real part. 

 The sum of the corresponding generalized eigenspaces are denoted $H,$ $E_-,$ and $E_+,$ so that
\be \label{space:dec} \C^N = H(t,x,\xi) + E_-(t,x,\xi) + E_+(t,x,\xi),
\ee
for all $t \leq T$ and $(x,\xi)$ in the support of $\psi.$ We denote $P_H,$ $P_{-},$ $P_{+}$ the associated eigenprojectors. By spectral separation, the projectors are pointwise as smooth as the solution, thus belong to $C^{[s - 1 - d/2], s -1 - d/2 - [s-1 -d/2]} S^0$ (see Appendix \ref{sec:symb} for notation on symbols and pseudo-differential operators). 
   The associated components of the unknown $v$ are 
\be \label{def:vEH}
v_H := \op(\psi^\flat P_H) v, \qquad v_{\pm} = \op(\psi^\flat P_{\pm}) v,
\ee
where we recall that $\psi^\flat$ is a cut-off with a smaller support than $\psi,$ so that $(1 - \psi)\psi^\flat \equiv 0.$ 
We denote
\be \label{def:MHE}
M_H := \psi^\sharp P_H \tilde A, \qquad M_\pm := \psi^\sharp P_\pm \tilde A.
\ee
The subsystems in $v_H,$ $v_-$ and $v_+$ are 
\be \label{eq:H+-}
 \d_t v_\star + \op(M_\star) v_\star = g_\star + R_\star v,
 \qquad \star \in \{H, -, +\},
 \ee
with notation
\be \label{def:gstar}
\ba 
g_\star & := \op(\psi^\flat P_\star) g, \\ 
R_\star & :=- \op(\psi^\flat \d_t P_\star) -  \op(\psi^\flat P_\star) \op(\psi^\sharp \tilde A) + \op(M_\star) \op(\psi^\flat P_\star).
\ea 
\ee   

\begin{lem} \label{lem:gstar} We have $g_\star \in L^\infty([0, T], H^{2 s - 2 - d/2}(\R^d)).$ The operator $R_\star$ maps $H^{s'}$ to itself, for any $s' \in \R,$ uniformly in $t \in [0,T].$
\end{lem}

The term $g_\star$ is akin to the ``out'' remainder terms in \cite{NT1}. 

\begin{proof} 
 By regularity of $\psi^\flat,$ and the posited regularity of $u,$ the operator $\op(\psi^\flat P_\star)$ maps $H^{s'}$ to itself, for any $s' \in \R$ (see Appendix \ref{sec:symb}). 
The map $g$ is defined by $g = f + f_1,$ which are defined in \eqref{def:f} and \eqref{def:f1}. In $f$ the first term is $\op(\psi) R^{para} \in H^{2 s - 2 -d/2},$ as seen in Section \ref{sec:paralin}. The second term in $f$ contributes to $g_\star$ the term
 \be \label{00} \op(\psi^\flat P_\star) [\op(\psi), \op(\tilde A)] u.\ee
 For the operator in \eqref{00}, we can use the composition result from Appendix \ref{sec:symb}. Since derivatives of $\psi$ are identically equal to 0 on the support of $\psi^\flat,$ we see that $\op(\psi^\flat P_\star) [\op(\psi), \op(\tilde A)]$ is linear bounded from $H^{s'}$ to $H^{s' - 1 + (s - 1 - d/2)},$ for any $s' \in \R.$ Thus the map in \eqref{00} belongs to $H^{2 s - 2 - d/2}.$  
The first term in $f_1$ contributes to $g_\star$ the term 
\be \label{01}
\op(\psi^\flat P_\star) \op(\tilde A) \big(\op(\psi^\sharp) \op(\psi) - \op(\psi)\big) u.
\ee  Here we may expand the composition $\op(\psi^\sharp) \op(\psi)$ up to any order, and find that $\op(\psi^\sharp) \op(\psi) - \op(\psi)$ is infinitely regularizing. Thus the map in \eqref{01} belongs to $H^{s'},$ for any $s'.$ 
The second term in $f_1$ contributes to $g_\star$ the term 
\be \label{02} \op(\psi^\flat P_\star)  (\op(\psi^\sharp \tilde A) - \op(\tilde A) \op(\psi^\sharp)\big) v.\ee
This term is seen to belong to $H^{2s - 2 - d/2}$ exactly like \eqref{00}. The bounds are uniform in time, since we assumed $\| u(t) \|_{H^s}$ to be bounded in time over $[0,T].$

 We turn to $R_\star.$ By the assumed regularity of $u,$ the symbol $\psi^\flat \d_t P_\star$ is bounded in $x.$ Indeed, by spectral separation, and regularity of $F,$ the symbol $P_\star$ is smooth in $(t,x,u,\d_x u),$ and, using \eqref{equation reference}:  
 $$ \d_t P_\star = \d_1 P_\star - \d_3 P \cdot F - \d_4 P \cdot \d_x F.$$
 Since $s > 2 + d/2,$ the map $\psi^\flat \d_x F,$ which involves two spatial derivatives of $u,$ belongs to $C^{[s - 2 - d/2], s - 2 - d/2 - [s-2-d/2]} S^0,$ and so does $\psi^\flat \d_t P_\star:$ 
 \be \label{reg:dtp} 
  \psi^\flat \d_t P_\star \in C^{[s - 2 - d/2], s - 2 - d/2 - [s-2-d/2]} S^0. 
 \ee 
 In particular, the norm $\| \psi^\flat \d_t P_\star \|_{0,0,[d/2] + 1}$ is finite, uniformly in time (here we used the symbolic norms introduced in Appendix \ref{sec:symb}), so that, by \eqref{action}, the norm $\| \op(\psi^\flat \d_t P_\star) \|_{H^{s'} \to H^{s'}}$ is bounded in time.   
 
 The other term in $R_\star$ is the difference 
$\op(\psi^\flat P_\star) \op(\psi^\sharp \tilde A) - \op(M_\star) \op(\psi^\flat P_\star).$ 
We observe that $\psi^\flat P_\star \psi^\sharp \tilde A - M_\star \psi^\flat P_\star = 0.$ Thus
\be \label{for:Rstar} \ba \op(\psi^\flat P_\star) \op(\psi^\sharp \tilde A) & - \op(M_\star) \op(\psi^\flat P_\star)
\\ & = (-i) \op\Big( \d_\xi (\psi^\flat P_\star) \d_x (\psi^\sharp \tilde A) - \d_\xi M_\star \d_x \psi^\flat P_\star\Big) \\ & + R_{s-1-d/2}(\psi^\flat P_\star,\psi^\sharp \tilde A) - R_{s-1-d/2}(M_\star, \psi^\flat P_\star), \ea\ee 
using notation for remainders introduced in Appendix \ref{sec:symb}. In the first line of the above right-hand side, we find a symbol that belongs to $C^{[s-2 - d/2], s - 2 - d/2 - [s-2-d/2]} S^0.$ Since $s - 2 > d/2,$ the corresponding operator is $H^{s'} \to H^{s'}$ bounded, uniformly in time over the compact $[0,T].$ In the second line, we find a remainder that enjoys the bound \eqref{remainder} with $m_1 + m_2 = 1,$ so that the remainder is order $1 - (s-1-d/2) < 0.$ In particular, is is $H^{s'} \to H^{s'}$ bounded. 
 \end{proof}

\subsection{The positive elliptic subsystem: forward-in-time propagation} \label{sec:ell:subsystem}

Consider the elliptic subsystem associated with eigenvalues with positive real parts, that is \eqref{eq:H+-} with $\star = +:$
\be \label{eq:+}
 \d_t v_+ + \op(M_+) v_+ =  g_+ + R_+ v, \qquad v_+(0) = \op(\psi^\flat P_+) \op(\psi) u_{in}.%
 \ee

The solution $v_+$ to \eqref{eq:+} is regularized as it is propagated forward in time:

\begin{prop} \label{prop:+}  
If the supports of $\psi$ and $\psi^\sharp$ are small enough and far enough from 0, and $T > 0$ is small enough, we have $v_+(t) \in C^0([0, T], H^{s + \rho}(\R^d)) \cap L^2([0,T], H^{s + \rho + 1/2}(\R^d)).$
\end{prop} 

Recall that $\rho$ is introduced in \eqref{def:rho}. In the above statement of Proposition \ref{prop:+}, the conditions on the supports mean $U_{x^0}^\sharp$ and $V_{\xi^0}^\sharp$ small enough, and the radius $r^\sharp$ large enough. See Section \ref{sec:loc} for the definition of $\Omega$ and $\Omega^\sharp.$ 

\begin{proof} 
{\it 1. Spectral bound.} By Assumption \ref{ass:ell}, choice of the initial datum $u_{in},$ and definition of the spectral projection $P_+,$ the spectrum of $M_+$ at $(x^0,\xi^0)$ has positive real part. By continuity, if the support of $\psi$ is chosen to be far away from 0, so that $\tilde A$ is a small perturbation of $A,$ we have a lower bound %
\be \label{low:+}
 \Re e \, M_+ \geq \g |\xi| >0, 
\ee
 in a whole conic (in $\xi$) region near $(x^0,\xi^0),$ for $t \in [0,T],$ if $T$ is small enough. We may assume that the support $\Omega^\sharp$ was chosen small enough so that \eqref{low:+} holds in the whole domain $\Omega^\sharp.$ 
 In \eqref{low:+}, $\Re e \, M_+ = (M_+ + M_+^\star)/2$ denotes the symmetric part of the matrix $M_+.$

\medskip

{\it 2. Regularization.} 
Let $j_\e$ be a regularizing kernel $j_\e: \R^d \to \R_+,$ such that 
 \be \label{reg:0} \| f - j_\e \star f \|_{L^\infty} \lesssim \e \| f \|_{W^{1,\infty}}, \qquad \mbox{for all $f \in W^{1,\infty}.$}\ee
 We let $M^\e_+ := M_+(t,x, j_\e \star u, j_\e \star \d_x u).$ 
 Since $s - 1 > 1 + d/2,$ this implies 
\be \label{reg:M} \langle \xi \rangle^{1 - |\b|} \| \d_\xi^\b (M_+ - M_+^\e) \|_{L^\infty} \lesssim \e.\ee
From \eqref{reg:M} we deduce first that the spectral bound \eqref{low:+} holds as well for $M_+^\e$ (up to an $O(\e)$ change in $\g_+$), and second (via \eqref{action}) the bound 
\be \label{approx:M} \| \op(M_+ - M_+^\e) \|_{H^{s' +1 } \to H^{s'}} \lesssim \e, \qquad \mbox{for all $s'.$}\ee
Spatial derivatives of $M_+^\e$ are stiff in $\e$ (below, we use symbolic norms $\| \cdot \|_{m,k,k'}$ defined in Appendix \ref{sec:symb}):
\be \label{stiff} \| M_+^\e \|_{1,k,k'} \lesssim \e^{-(k-1)}.\ee 
Above, the exponent is $k-1$ and not $k$ since the argument $(u,\d_x u)$ of $M_+^\e$ belongs to $W^{1,\infty}.$ 

\medskip

{\it 3. The $L^2$ estimate.} 
  In a next step, we apply $J D_x^{s + \rho}$ to the equation then perform an $L^2$ estimate. Here $J$ is the regularizing operator $J = (1 + \delta |D|^2)^{-(1 + \rho)/2} \in S^{- (1 + \rho)},$ which depends on $\delta > 0.$ We let $w = J D^{s + \rho} v_+ \in H^{1}$ and observe that $w$ satisfies
  \be \label{eq:w}
  \d_t w + \op(M_+) w = J D^{s+\rho} (g_+ + R_+ v) + \Gamma v_+, \qquad \Gamma := [\op(M_+),  J D^{s+\rho}].\ee 
In particular, $\d_t w \in L^2,$ and we have %
\be \label{L2:est} \ba \frac{1}{2} \Re e \, \d_t (\| w \|_{L^2}^2) + \g \Re e \, ( \op(|\xi|\psi^\sharp) w,  w)_{L^2} & + \Re e \, \big( \op(M_+^\e - \g|\xi| \psi^\sharp) w, w \big)_{L^2} = \Re e \, f,\ea\ee
with
$$ \ba f & = \big( \op(M_+^\e - M_+) w, w \big)_{L^2} + \big(J D_x^{s+\rho} (g_+ + R_+ v), w )_{L^2} + (\Gamma v_+, w)_{L^2} =: \sum_{1 \leq j \leq 3} f_j.\ea$$
Consider first the left-hand side of \eqref{L2:est}. Since $w$ is defined in terms of $v_+,$ which is defined in terms of $\psi^\flat,$ the truncation $\psi^\sharp$ does not contribute much to $\op(|\xi| \psi^\sharp) w.$ That is, 
$$  \op(|\xi| \psi^\sharp) w = |D| w + \tilde w,$$
with
$$ 
 \ba \tilde w := R_1(\psi^\sharp J |\xi|^{1 + s + \rho}, \psi^\flat P_+) v - R_1(|\xi|, J |\xi|^{s + \rho} \psi^\flat P_+) v - |D| R_1(\psi^\sharp J |\xi|^{s + \rho}, \psi^\flat P_+) v,
 \ea
 $$ 
 and for fixed $\delta > 0,$ we have $\tilde w \in L^\infty H^1.$ %
Thus the left hand side of \eqref{L2:est} is 
\be \label{L2:est:1} \frac12 \Re e \, \d_t (\| w \|_{L^2}^2) + \g \| |D|^{1/2}  w \|_{L^2}^2 + \g \Re e \, (\tilde w, w )_{L^2} + \Re e \, \big( \op(M_+^\e - \g |\xi| \psi) w, w \big)_{L^2}.\ee

\medskip

{\it 4. G\r{a}rding's inequality.} We bound the contribution of $M_+^\e$ to \eqref{L2:est:1} with G\r{a}rding's inequality. First we compare $\op(M_+^\e)$ to the associated pseudo-differential operator:
\be \label{op:pdo} \big\| (\op(M_+^\e) - \pdo(M_+^\e)) w \big\|_{L^2} \lesssim C(M_+) \| w \|_{L^2}.\ee
The \eqref{op:pdo} estimate follows from Remark A.4 (and Proposition A.3 and its proof) in \cite{NT1}, based on the results of \cite{Lannes}. Above, $C(M_+)$ is 
$$ \ba C(M_+) & = \sup_{ \begin{smallmatrix} |\b| \leq 2 [d/2] + 2 \\ \xi \in \R^d \end{smallmatrix} } \langle \xi \rangle^{|\b| - 1} \| \d_\xi^\b ((1 - \tilde \phi_0(D_x)) M_+) \|_{H^{1 + d/2}} \\ &\qquad + \| (1 - \tilde \phi_0(D_x)) \phi_0(2^{- N_0} D_x) M_+(\cdot,\xi) \|_{1, 1 + [d/2], 1 + [d/2]},\ea$$
where $\tilde \phi_0$ and $\phi_0$ are low-frequency cut-offs, and $N_0 \geq 3.$ Since $(u, \d_x u) \in H^{s-1},$ with $s - 1 > 1 + d/2,$ the first term in $C_+$ is finite, uniformly in $t:$
$$ \langle \xi \rangle^{|\b| - 1} \| \d_\xi^\b ((1 - \tilde \phi_0(D_x)) M_+) \|_{H^{1 + d/2}}  \lesssim \sup_{t \in [0,T]} C(\| (u, \d_x u)(t)\|_{L^\infty}) \| (u,\d_x u)(t)\|_{H^{1 + d/2}} < \infty.$$
 The second term in $C(M_+)$ is a sup norm involving $(x,\xi)$ derivatives of $M_+$ (see the definition of symbolic norms in Appendix \ref{sec:symb}). This term is finite too, since the inverse Fourier transform of the smooth frequency cut-off $(1 - \tilde \phi_0(\xi)) \phi_0(2^{- N_0} \xi)$ belongs to the Schwartz class, so that its derivatives belong to $L^1:$ 
 $$ \| (1 - \tilde \phi_0(D_x)) \phi_0(2^{- N_0} D_x) M_+(\cdot,\xi) \|_{1, 1 + [d/2], 1 + [d/2]} \lesssim \sup_{t \in [0,T]} C(\| (u, \d_x u)(t)\|_{L^\infty}) < \infty.$$
 By the spectral lower bound \eqref{low:+} and G\r{a}rding's inequality: 
$$  \Re e \, \big( \pdo(M_+^\e - \g |\xi| \psi^\sharp) w,w\big)_{L^2} + \| M_+^\e\|_{1,d_\star,d_\star} \| w \|_{L^2}^2 \geq 0,$$
where $d_\star$ is an integer which depends on the space dimension. (The proof of Theorem 1.1.26 in \cite{Lerner} shows that it is possible to take $d_\star = 3 + d.$ We do not need a precise value for $d_\star$ here, as we will simply choose $\e$ to be small enough later on.) By \eqref{stiff}, this implies 
$$  \Re e \, \big( \pdo(M_+^\e - \g_+ |\xi| \psi^\sharp) w,w\big)_{L^2} + C \e^{-(d_\star -1)} \| w \|_{L^2}^2 \geq 0,$$
for some $C > 0$ which does not depend on $\e.$ All bounds in this fourth step do not involve $J,$ hence do not depend on $\delta.$ 

\medskip

{\it 5. Remainder bounds.} There are seven remainder terms: the four terms in $f,$ defined in step 3, the term in $\tilde w$ also defined in step 3, and the G\r{a}rding remainder and the para- to pseudo- remainder from step 4. 

Only $w$ is involved in the terms $f_1,$ $f_2,$ $f_3,$ and in both remainders from step 4. By \eqref{approx:M}, we have $|f_1| \lesssim \e \| w\|_{H^{1/2}}^2.$ By Lemma \ref{lem:gstar}, we have 
$$ |f_2| \leq \| D^{-\rho} J D^{s + \rho} (g_+ + R_+ v) \|_{L^2} \| D^\rho w \|_{L^2} \lesssim \| v \|_{H^s} \| D^\rho w \|_{L^2},$$
 uniformly in $\delta.$ Taken together, the remainders from step 4 contribute at most $(C(M_+) + C \e^{-(d_\star - 1)}) \| w \|_{L^2}^2,$ where $C > 0$ does not depend on $\e,\delta.$

We turn to $f_3,$ defined in step 3. We have 
$$ \Re  e \, \big( [ \op(M_+), J D^{s + \rho} ] v_+, w \big)_{L^2} \leq \| D^{-\rho} [\op(M_+), J D^{s + \rho}] v_+ \|_{L^2}  \| D^\rho w\|_{L^2}.$$
 By \eqref{composition:para}, 
 $$ D^{-\rho} [\op(M_+), J D^{s + \rho}] = D^{-\rho} R_1(M_+, J |\xi|^{s + \rho}) - D^{-\rho} R_1(J |\xi|^{s + \rho}, M_+)),$$
 and, by \eqref{remainder},
 $$ \| D^{-\rho} R_1(M_+, J |\xi|^{s + \rho}) v_+ \|_{L^2} \leq \| R_1(M_+, J |\xi|^{s + \rho}) v_+ \|_{H^{-\rho}} \lesssim \| v_+ \|_{H^s},$$
 uniformly in $\delta,$ since $J$ and its (frequency) derivatives are bounded uniformly in $\delta,\xi.$ The other remainder $R_1(\dots,\dots)$ above is handled identically. Thus
 $$ |f_3| \lesssim \| D^\rho w\|_{L^2} \| v_+\|_{H^s},$$
 uniformly in $\e$ and $\delta.$  

The dot product $(\tilde w, w)_{L^2}$ remains. All three terms in $\tilde w$ are similar, so that we will focus on the first. It is handled just like the commutator involving $M_+$ in $f_4:$ 
$$ \big| \big( R_1(\psi^\sharp J |\xi|^{1 + s + \rho}, \psi^\flat P_+) v, w \big)_{L^2} \big| \leq \| D^{-\rho} R_1(\psi^\sharp J |\xi|^{1 + s + \rho}, \psi^\flat P_+) v \|_{L^2} \| D^\rho w \|_{L^2},$$
and, by \eqref{remainder}, 
$$  \| D^{-\rho} R_1(\psi^\sharp J |\xi|^{1 + s + \rho}, \psi^\flat P_+) v \|_{L^2} \lesssim \| v \|_{H^s},$$
uniformly in $\delta.$ Thus
$$ | (\tilde w, w )_{L^2}| \lesssim \| v \|_{H^s} \| D^\rho w \|_{L^2}.$$
  
 \medskip

{\it 6. Bound in $L^2 H^{s + \rho + 1/2}.$} Summing up, we obtain
\be \label{est:3} \ba \frac12 \Re e \, \d_t (\| w \|_{L^2}^2) & + \g \| |D|^{1/2}  w \|_{L^2}^2  \\ &  \leq C \Big( \e \| w \|_{H^{1/2}}^2 + \| v \|_{H^s} \| w \|_{H^{\rho}} + (1 + \e^{-(d_\star -1)}) \| w \|_{L^2}^2 + \| w \|_{L^2} \Big). \ea\ee 
Up to a constant, the terms in the above right-hand side can be absorbed into the $\g \| |D|^{1/2} w \|_{L^2}^2$ term in the left-hand side. Indeed, we have, using $\rho \leq 1/2,$ 
$$ \e \| w \|_{H^{1/2}}^2 + \| v \|_{H^s} \| w \|_{H^{\rho}} \leq (\e + \zeta) \| w \|_{H^{1/2}}^2 + \frac{1}{4 \zeta} \| v \|_{H^s}^2,$$
and 
$$ \| w \|_{H^{1/2}}^2 \leq 2 \| v_+\|_{L^2}^2 + \| |D|^{1/2} w \|_{L^2}^2,$$
so that 
$$ \e \| w \|_{H^{1/2}}^2 + \| v \|_{H^s} \| w \|_{H^{\rho}} \leq (2\e + 2 \zeta + \frac{1}{4\zeta}) \| v_+ \|_{H^s}^2 + (\e + \zeta) \| |D|^{1/2} w \|_{L^2}^2.$$ 
 We now choose $\e$ and $\zeta$ so that $C(\e + \zeta) \leq \g/2,$ where $C$ is the constant that appears in the right-hand side of \eqref{est:3}. This constant depends only on $F,$ the datum, and the spatial dimension $d.$ Then, the estimate \eqref{est:3} becomes 
\be \label{est:4} \ba \frac12 \Re e \, \d_t (\| w \|_{L^2}^2) + \frac{\g}{2} \| |D|^{1/2}  w \|_{L^2}^2 \leq C (1 + \e^{-(d_\star -1)}) \| w \|_{L^2}^2 \big) + C_{\e,\zeta} (1 + \| v_+\|_{H^s}^2),\ea\ee
for some $C_{\e,\zeta} > 0.$ Since the $H^s$ norm of $v_+$ is bounded over $[0,T],$ this implies, by Gronwall,
 $$ \mbox{for all $t \in [0,T],$} \quad \| w(t)\|_{L^2}^2 \lesssim C(\e,T, \| v_+\|_{L^2([0,T], H^s)}^2),$$
Using \eqref{est:4} again, and denoting $w_\delta = w,$ we see now that the sequence $(w_{\delta})_{\delta > 0}$ is bounded in the Banach space $L^2([0,T],H^{1/2}).$ Up to a subsequence, $w_\delta$ thus converges weakly in $L^2([0,T], H^{1/2}).$ Besides, $w_\delta$ converges in the sense of distributions to $D^{s + \rho} v_+.$ This proves $v_+ \in L^2([0, T], H^{s + \rho + 1/2}).$ 

\medskip

{\it 7. Conclusion.} From $v_+ \in L^2 H^{s + \rho + 1/2}$ and the equation in $v_+,$ we deduce that $\d_t v_+ \in L^2 H^{s + \rho - 1/2}.$ Thus $v_+ \in C^0 H^{s + \rho},$ which concludes the proof. 
\end{proof}

\subsection{The negative elliptic subsystem: backward-in-time propagation} \label{sec:bwd}

Consider the elliptic subsystem associated with eigenvalues with negative real parts, that is \eqref{eq:H+-} with $\star = -:$
\be \label{eq:-}
 \d_t v_- + \op(M_-) v_- =  g_- + R_- v, \qquad v_-(0) = \op(\psi^\flat P_-) \op(\psi) u_{in}.%
 \ee

The unknown $v_-$ is regularized as it is propagated backward in time:

\begin{prop} \label{prop:-}  
Under the assumptions of Proposition {\rm \ref{prop:+}}, given $0 < T_- < T,$ we have $v_-(t) \in C^0([0, T_1], H^{s + \min(1/2,s-2-d/2)}(\R^d)) \cap L^2([0,T_1], H^{s + \min(1/2,s-2-d/2) + 1/2}(\R^d)).$
\end{prop} 

\begin{proof} Let $T_1$ such that $0 < T_- < T_1 < T,$ and 
$$  v_{bwd}(t) = v_-(T_1 - t), \qquad t \in [0, T_1],$$
so that $v_{bwd} \in C^0([0,T_1], H^s(\R^d))$ solves
\be \label{eq:w-} \left\{\ba
   \d_t v_{bwd} & = \op(M_-(T_1 - t)) v_{bwd} - g_-(T_1 - t) - R_-(T_1 - t) v(T_1 - t), \\ v_{bwd}(0) & = v(T_1).\ea \right.\ee
By Assumption \ref{ass:ell}, choice of the initial datum $u_{in},$ and definition of the spectral projection $P_-,$ the spectrum of $M_-$ at $(x^0,\xi^0)$ has negative real part. By reality of $F$ and the datum, the spectrum of $A$ is symmetric with respect to the imaginary axis, so that the real part of the spectrum of $M_-$ is the opposite of the real part of the spectrum of $M_+.$ Thus under the conditions of Proposition \ref{prop:+}, we have
\be \label{low:-}
 \Re e \, M_-(T_1 - t) \leq - \g |\xi| < 0, \qquad \mbox{for $(x,\xi) \in \Omega^\sharp,$ and $t \in [0,T_2].$}
\ee
We may now repeat the arguments of the proof of Proposition \ref{prop:+}. The only significant change is the fact that $M_-$ is now on the right-hand side of the equation. The sign change for $g_-$ and $R_-,$ and the evaluation at $T_2 - t$ rather than $t$ do not matter: indeed, in the course of Proposition \ref{prop:+} these terms are bounded uniformly in time. 

From \eqref{low:-}, we deduce, by G\r{a}rding's inequality,
$$ \Re e \, \big( \pdo \big(M_-^\e(T_1 -t) + \g |\xi| \psi\big) z, z \big)_{L^2} \leq C(\e) \| z \|_{L^2}^2,$$
  for all $z \in L^2(\R^d),$ where $M_-^\e$ is a regularized symbol (see the proof of Proposition \ref{prop:+}).

  We now differentiate the equation in $v_{bwd}$ and perform an $L^2$ estimate. Here a difference with Proposition \ref{prop:+} appears: the datum for $v_{bwd}$ belongs to $H^s,$ not $H^\s.$ Thus we let $w_- := J_- D^{s} v_{bwd},$ with $J_- := (1 + \delta |D|^2)^{-1/2},$ and for $w_-$ we find, with the above G\r{a}rding inequality, the $L^2$ estimate  
 $$ \frac{1}{2} \d_t (\| w_- \|_{L^2}^2) + \g \| |D|^{1/2} w_-\|_{L^2}^2 \leq C(\e) \| w_-\|_{L^2}^2 + \mbox{other error terms},$$
 and the ``other error terms", above, are bounded in exactly the same way as the error terms in the proof of Proposition \ref{prop:+}. This implies a bound analogous to \eqref{est:4}:
$$  \ba \frac12 \Re e \, \d_t (\| w_- \|_{L^2}^2) + \frac{\g}{2} \| |D|^{1/2}  w_- \|_{L^2}^2 \leq C (1 + \e^{-(d_\star -1)}) \| w_- \|_{L^2}^2 \big) + C_{\e} (1 + \| v \|_{H^s}^2),\ea
$$ %
 where the constant $C$ is independent of $\e$ and $\delta.$ We conclude as in the proof of Proposition \ref{prop:+} that $v^{bwd}_- \in C^0([0,T_1], H^{s}) \cap L^2([0,T_1],H^{s + 1/2}).$ (Here we used the arguments of Proposition \ref{prop:+}, with
 $\rho_- := \min(s - s, s-2-d/2, 1/2) = 0,$
 since the datum has regularity $H^s.$)

 Thus for almost all $t \in [0,T_1],$ we have $v_{bwd}(t) \in H^{s + 1/2}.$ Let $T_2 \in (0, T_1 - T_-)$  be such a time. We may repeat the above argument with the initial-value problem starting at $T_2$ with $v_{bdw},$ and find that $v_{bwd} \in C^0([0, T_2], H^{s + \min(s - 2 - d/2,1/2)}) \cap L^2([0,T_2], H^{s + \min(s-2-d/2,1/2) + 1/2}).$ (Here, we used the arguments of Proposition \ref{prop:+}, with
 $\tilde \rho_- = \min(s + 1/2 - s, s-  2 - d/2, 1/2) = \min(s- 2 -d/2,1/2),$
 since the datum has regularity $s + 1/2.$) By definition of $v_{bwd}$ at the beginning of this proof, this implies the same regularity for $v_-$ over the time interval $[0, T_-],$ and the result. 

\end{proof}

\subsection{The hyperbolic subsystem}  \label{sec:hyp}

Consider the hyperbolic subsystem associated with purely imaginary eigenvalues, that is \eqref{eq:H+-} with $\star = H:$
\be \label{eq:H}
 \d_t v_H + \op(M_H) v_H =  g_H + R_H v, \qquad v_H(0) = \op(\psi^\flat P_H) \op(\psi) u_{in}.%
 \ee 
 
By Assumption \ref{ass:ell}, the purely imaginary eigenvalues of $A$ at $(0,\o^0)$ are semi-simple and have  constant multiplicity. This implies existence of a symmetrizer for \eqref{eq:+}. As a consequence, in short time the regularity of the solution $v_H$ is the minimum of the regularities of the datum and source:

\begin{prop} \label{prop:H} Under the conditions of Proposition {\rm \ref{prop:+},} for some $0 < T_H \leq T,$ we have $v_H \in C^0([0,T_H], H^{s + \rho}(\R^d)),$ where $\rho$ is defined in \eqref{def:rho}. 
\end{prop}

\begin{proof} {\it 1. Spectral decomposition.} The eigenvalues of $M_H$ at $(0,\o^0)$ are semi-simple and have constant multiplicity. As a consequence, the branches of eigenvalues $i \l_j \in i \R $ of $M_H$ and the corresponding eigenprojectors $P_j$ are locally as regular as $(u,\d_x u),$ and we have 
$$ M_H(t,x,\xi) = \sum_{1 \leq j \leq j_0} i \l_j(t,x,\xi) P_j(t,x,\xi),$$
for some $j_0$ equal to the number of distinct eigenvalues of $M_H$ at $(0,\o^0).$ 
We may assume that the neighborhoods $U_{x^0}^\sharp$ and $V_{\xi^0}^\sharp$ of $x^0$ and $\xi^0$ are small enough and $r^\sharp$ is large enough so that the above decomposition holds in $[0, T_H] \times \Omega^\sharp,$ for some $0 < T_H \leq T.$ 

\medskip

{\it 2. A linear, para-differential, hyperbolic system with $W^{1,\infty}$ coefficients.} We let
\be \label{def:fH}
 f_H := g_H + R_H(v_- + v_+) \in L^\infty([0, T_+], H^{s + \rho}),
\ee
by Lemma \ref{lem:gstar} and Propositions \ref{prop:+} and \ref{prop:-}. 
We define the linear operator 
\be \label{def:L_H}
L_H := \d_t + \op(M_H) - R_H.
\ee
 The hyperbolic initial-value problem \eqref{eq:H} takes the form 
\be \label{eq:LH} L_H v_H = f_H, \qquad v_H(0) \in H^{\s}.\ee
The key is that the coefficients of $M_H,$ which is order one, belong to $W^{1,\infty}$ by regularity of $u.$ 

\medskip

{\it 3. Sobolev estimate.} For all $f \in L^2,$ we have 
$$ \Re e \, \big( L_H f, f \big)_{L^2}  = \frac{1}{2} \d_t \| f \|_{L^2}^2 + \sum_{1 \leq j \leq j_0} \Re e \, \big( \op(i\l_j P_j) f, f)_{L^2} - \Re e \, \big( R_H f ,f \big)_{L^2},$$
and, by reality of the $\l_j,$ 
$$ \ba \Re e \, \big( \op(i\l_j P_j) f, f)_{L^2} & = \frac{1}{2} \big(\, \big( \op(i \l_j P_j) + \op(i \l_j P_j)^\star) \big) f, f \, \big)_{L^2} = \frac{1}{2} (R_1^\star(i \l_j P_j) f, f)_{L^2}, \ea $$
using notation $R_1^\star(\dots)$ for the remainder in the description of an adjoint operator (see Appendix \ref{sec:symb}). By regularity of the spectral projectors $P_j$ and \eqref{remainder:adjoint}, the $L^2 \to L^2$ norm of $R_1^\star(i \l_j P_j)$ is bounded, uniformly in $t.$ Here we used the fact that $\l_j$ and $P_j$ depend smoothly on $(u, \d_x u)$ which belong to $W^{1,\infty}$ in a neighborhood of $x^0.$ Thus 
$$ \big|\, \Re e \, \big( \op(i\l_j P_j) f, f)_{L^2}\, \big| \lesssim \|f \|_{L^2}^2,$$
uniformly in $t \in [0, T_H],$ 
so that, using Lemma \ref{lem:gstar} for a bound of $R_H,$ 
\be \label{est:H}
\frac{1}{2} \d_t (\| f \|_{L^2}^2 \leq \Re e \, (L_H f, f)_{L^2} + C \| f \|_{L^2}^2, 
\ee  
for all $f \in L^2$ and some constant $C > 0$ which does not depend on $t.$ By the $W^{1,\infty}$ regularity of the coefficients of $M_H,$ and the composition result \eqref{composition:para}-\eqref{remainder}, we have 
\be \label{commut:H}
 \| [D^{s + \rho}, \op(M_H)] f \|_{L^2} \lesssim \| f \|_{H^{s + \rho}}, \qquad \mbox{for all $f \in H^{s +\rho}.$}
\ee
 Consider now the commutator $[D^{s + \rho}, R_H].$ Here we use the proof of Lemma \ref{lem:gstar}. The first term in $R_H$ involves $\d_t P_H,$ which is less regular than $M_H:$ its regularity is given in \eqref{reg:dtp}. Thus, by \eqref{composition:para} 
 $$  [\op(\psi^\flat \d_t P_H), D^{s + \rho}] = - R_{s-2-d/2}(D^{s + \rho}, \psi^\flat \d_t P_H),$$
 so that, by \eqref{remainder},
 \be \label{commut:J} \| [\op(\psi^\flat \d_t P_H), D^{s + \rho}] v_H \|_{L^2} \lesssim \| v_H \|_{H^{s + \rho - s - 2 - d/2}},\ee
 and $\rho \leq s - 2 - d/2,$ by definition of $\rho$ \eqref{def:rho}. Consider now the other term in $R_H.$ The symbol in the first line of the right-hand side of \eqref{for:Rstar} belongs to the same class as $\psi^\sharp \d_t P_H$ \eqref{reg:dtp}, hence its commutator with $J D^{s + \rho}$ is bounded just like in \eqref{commut:J}. The term in the second line of the right-hand side of \eqref{for:Rstar} is linear and bounded from $H^{s'}$ to $H^{s' + 1 - (s - 1 -d/2)},$ by \eqref{remainder}. Since $s' + 1 - (s - 1 - d/2) < 0,$ by assumption on $s,$ this implies that the commutator of that term with $D^{s + \rho}$ is linear bounded from $H^{s + \rho}$ to $L^2.$    We conclude that
 \be \label{commut:R+}
  \| [  D^{s + \rho}, R_H ] f \|_{L^2 }\lesssim \| f \|_{H^{s + \rho - (s - 2 - d/2)}}, \quad \mbox{for all $f \in H^{s + \rho - (s - 2 - d/2)},$} \ee
  with an implicit constant which is independent of $t.$

From there, we deduce the $H^{s + \rho}$ estimate
\be \label{est:H:s'}
\frac{1}{2} \d_t (\| f \|_{H^{s + \rho}}^2) \leq \Re e \, (L_H f, f)_{H^{s + \rho}} + C \| f \|_{H^{s+ \rho}}^2. 
\ee 

\medskip

{\it 4. Existence and uniqueness of a solution in $C^0([0,T_H], H^{s'}).$} From \eqref{est:H:s'} and the regularity of $f_H$ \eqref{def:fH}, it is classical to deduce the existence and uniqueness of a solution to the initial-value problem \eqref{eq:LH} which is continuous in time over the whole interval $[0,T_H]$ (recall, this initial-value problem is linear since the existence of $u$ is posited in the first place) with values in $H^{s + \rho}.$ 

By uniqueness, that solution coincides with $v_H.$ Thus $v_H$ has the announced regularity. 
\end{proof} 

\subsection{Endgame}  \label{sec:endgame}

By \eqref{space:dec} and \eqref{def:vEH}, we have $v = v_+ + v_- + v_H.$ By Propositions \ref{prop:+}, \ref{prop:-} and \ref{prop:H}, this proves that $v \in C^0([0, \min(T_-,T_H)], H^{s + \rho}(\R^d)).$ 

As noted at the beginning of Section \ref{sec:loc}, all frequencies $\xi \in \S^{d-1}$ are elliptic. So associated with every $\xi$ we have $U(\xi) \subset \R^d$ a neighborhood of $x^0$ and $V(\xi) \subset \S^{d-1}$ a neighborhood of $\xi$ so that the above analysis holds, meaning that we have $v_\xi := \op(\psi^\flat_\xi) u \in C^0([0, T_\xi], H^{s + \rho}(\R^d)),$ for some $T(\xi) > 0,$ where $\psi^\flat_\xi$ is an space-frequency cut-off, defined in terms of $U(\xi)$ and $V(\xi)$ just like in Section \ref{sec:loc}. 

By compactness of the sphere, a finite number of the associated $V(\xi)^\flat$ cover $\S^{d-1}.$ Let $V_1^\flat, \dots, V_n^\flat$ be this finite open cover. Let $U = \cap_{1 \leq i \leq n} U_i^\flat:$ the intersection of the associated neighhorhoods of $x^0,$ itself an open neighborhood of $x^0.$ The associated space-frequency cut-offs are denoted $\psi_i, \psi_i^\flat, \psi_i^\sharp.$ Let $T_0 = \min_i T_i,$ where $T_i$ is such that $v_i := \op(\psi_i^\flat) u \in C^0([0,T_i], H^{s  +\rho}).$ 

From there we argue as follows in order to prove that $u$ itself, in restriction to a small neighborhood of $x^0,$ belongs to $C^0([0, T_0], H^{s + \rho}).$ By Remark A.4 from \cite{LNT} (based on the analysis of \cite{Lannes}), for any symbol $a$ of order $m$ such that the upper bound below is finite, we have
$$ \ba \| \op(a) - \pdo(a) \|_{H^{m + d/2} \to H^{s'}} & \lesssim \sup_{ \begin{smallmatrix} |\b| \leq 2 [d/2] + 2 \\ \xi \in \R^d \end{smallmatrix} } \langle \xi \rangle^{|\b| - m} \| \d_\xi^\b ((1 - \tilde \phi_0(D_x)) a) \|_{H^{s'}} \\ & \quad + \| (1 - \tilde \phi_0(D_x)) \phi_0(2^{- N_0} D_x) a(\cdot,\xi) \|_{m, 1 + [d/2], 1 + [d/2]},\ea $$
for any $s',$ where $\phi_0$ and $\tilde \phi_0$ are smooth low-frequency cut-offs (compactly supported in a neighborhood of the zero frequency). Applied to $\psi_i^\flat,$ this gives the bound 
$$ \| \op(\psi_i^\flat) - \pdo(\psi_i^\flat) \|_{H^{d/2} \to H^{s'}} < \infty.$$
Thus 
$$ (\op(\psi_i^\flat) - \pdo(\psi_i^\flat) ) u \in C^0([0, T_0], H^{s + \rho}),$$
In particular,
$$ \sum_i \pdo(\psi_i^\flat) u \in C^0([0, T_0], H^{s + \rho}).$$
Now we may choose each $\psi_i^\flat$ to be a tensor product: $\psi_i^\flat  = \psi_{i1}(x) \psi_{i2}(\xi).$ For each $i,$ the spatial cut-off $\psi_{i1}$ is identically equal to 1 in a neighborhood of $x^0.$ It is not restrictive to assume that $U$ is small enough so that $\psi_{i1} \equiv 1$ on $U.$ Thus we have 
\be \label{for:endgame} \Big( \sum_i \psi_{i2}(D_x) u \Big)_{|U} \in C^0([0, T], H^{s + \rho}).\ee 
It is not restrictive to assume that $V_i^\flat$ is small enough so that $\psi_{i2} \equiv 1$ on $V^\flat_i.$ Thus $\chi := \sum_i \psi_{i2}$ is a high-frequency cut-off, such that
\be \label{chi} c_1 \leq \chi(\xi) \leq c_2, \quad \mbox{for all $|\xi| \geq r,$ for some $r > 0.$}\ee
Let $\chi_{LF}$ be compactly supported and such that $\chi_{LF} \equiv 1$ in a neighborhood of $|\xi| \leq r.$ Since $u \in C^0 H^{s},$ and $\chi_{LF}$ is compactly supported, we have $\chi_{LF}(D) u \in C^0 H^{s + \rho}.$ Besides, by \eqref{chi} and definition of $\chi_{LF},$ the Fourier multiplier $(1 - \chi_{LF}) \chi^{-1}(D)$ is well-defined and maps $H^{s'}$ to itself, for any $s'.$ Thus  
 $$ u_{|U} = (\chi_{LF} u)_{|U} + \Big( (1 - \chi_{LF}(D)) \chi^{-1}(D) \chi(D) u\Big)_{|U} \in C^0([0,T_0], H^{s +\rho}).$$ 
We may now repeat the argument, and after a few iterations reach a contradiction with Theorem 1 from \cite{NT1}, as explained in the first paragraph of this proof. 

\section{Proof of Theorem \ref{th:transition}: weak defects of hyperbolicity} \label{sec:weak}

 The main difference with the proof of Theorem \ref{th:main} is seen in the elliptic equations: the real parts of the eigenvalues are $O(t),$ hence the elliptic operator is $O(t),$ meaning a weaker regularization in forward time for the positive elliptic equations. The regularization index is here 
 \be \label{def:rho:bif} 
 \nu < 1/2,
\ee
meaning that $1/2 - \nu$ is positive and arbitrarily small. 

\subsection{The branching eigenvalues and eigenvectors} \label{sec:spec}

 We sum up here the spectral situation that Assumption \ref{ass:bif} entails. Details are given in Section 4.1 in \cite{NT1}.

 We have real eigenvalues that are simple and stay real, locally in $(t,x,\xi).$ The regularity of these eigenvalues in $(t,x,\xi)$ equals the one of the principal symbol $A$ evaluated at $(u,\d_x u).$ We denote $\underline A(t,x,\xi) = A(t,x,u(t,x), \d_x u(t,x), \xi).$ The symbol $\underline A$ is $C^1$ in $t,$ $C^{2,1/4}$ in $x$ (via the system \eqref{equation reference} and the condition on $s$ in Theorem \ref{th:transition}), and smooth in $\xi.$ 
 
 Besides, we have pairs of branching eigenvalues which leave the real axis at $t = 0.$ 
 
 At $t = 0,$ the regularity in $(x,\xi)$ of these eigenvalues is equal to the regularity of $\underline A_{|t = 0}$ evaluated at $(u_{in},\d_x u_{in}),$ since they have constant multiplicity as $(x,\xi)$ is varied. 
 
 For $t > 0,$ the branching eigenvalues are simple, hence their regularity is equal to the regularity of $\underline A.$
 
  The branching eigenvalues $\mu$ are time-differentiable at $t = 0,$ with $\d_t \mu \neq 0.$ Besides, $\d_t \mu$ is continuous in $(t,x)$ and smooth in $\xi.$ 
  
  To the branching eigenvalues are associated eigenvectors which have the same regularity as the branching eigenvalues. 

\subsection{Localization} \label{sec:loc:bif}  

 We use notation from Sections \ref{sec:paralin} and \ref{sec:loc}. Defining $v = \op(\psi) u$ as in \eqref{def:v}, we arrive at 
 \be \label{eq:v:bif} \left\{\begin{aligned}
 \d_t v + \op(i A) v & =  \op(\psi) R^{para} - [\op(\psi), \op(i A)] u - \op(\psi) \op(\d_3 F) u, \\ v(0)  & = \op(\psi) u_{in}. \end{aligned}\right.
 \ee 
 The only difference between \eqref{eq:v:bif} and \eqref{ivp:v} is that in \eqref{eq:v:bif}, the term $\d_3 F$ is seen as being part of the small ``source'' in the right-hand side, not the propagator. We may write
 $$ \op(\psi) \op(\d_3 F) u = \op(\d_3 F) v + [\op(\psi), \op(\d_3 F)] u,$$
 so that, defining the ``source'' term $f$ as in \eqref{def:f}, 
the initial-value problem \eqref{eq:v:bif} takes the form
\be \label{eq:v:2} \begin{aligned}
 \d_t v + \op(i A) v = f - \op(\d_3 F) v, \quad v(0) = \op(\psi) u_{in}. \end{aligned}\ee
Just like in Section \ref{sec:loc}, we may introduce a space-frequency truncation in the symbol of the propagator in \eqref{eq:v:2}, leading to 
\be \label{eq:v:3} \d_t v + \op(i \psi^\sharp A) v  = g - \op(\d_3 F) v, \quad v(0) = \op(\psi) u_{in},\ee
where $g = f + f_1,$ with $f_1$ defined in \eqref{def:f1} where $\tilde A$ is replaced by $i A.$

\subsection{Diagonalization} \label{sec:diag}

 At this point, unlike in the proof of Theorem \ref{th:main}, but as in the proof of Theorem 2 in \cite{NT1}, we use $Q$ a diagonalization matrix for $A.$ The regularity of $Q$ in $(t,x,\xi)$ is described in Section \ref{sec:spec}: the symbol $Q$ is as regular as $\underline A$ in $(t,x,\xi),$ and the symbol $\d_t Q$ is continuous in $(t,x)$ and smooth in $\xi.$ We extend the locally defined symbol $Q$ into a globally invertible symbol as in Appendix C of \cite{LNT}.  
We let 
\be \label{def:w} v_Q = \op(\psi^\flat Q) v,
\ee and find that the initial-value problem in $v_Q$ is 
\be \label{ivp:w} \d_t v_Q + \op(i \psi^\sharp Q A Q^{-1}) v_Q = g_Q + R_Q v, \qquad v_Q(0)  = \op(\psi^\flat Q) v(0). \ee 
where %
\be \label{def:gQ}
g_Q := \op(\psi^\flat Q) g, \qquad R_Q v := - \op(\psi^\flat Q) \op(\d_3 F) v - \op(\psi^\flat \d_t Q) v + R v,
\ee 
with notation
\be \label{def:RQ}
\ba R & = \big( R_1(i \psi^\flat Q A, Q^{-1}) + R_1(\psi^\flat Q, i \psi^\sharp A ) \op(Q^{-1}) \big) \op(Q) \\ & + \op(i \psi^\sharp Q A Q^{-1}) R_1(\psi^\flat, Q) + R_{1}(i \psi^\sharp Q A Q^{-1}, \psi^\flat) \op(Q) \\ & + \op(\psi^\flat Q) \op(i \psi^\sharp A) (\op(Q^{-1}) \op(Q) - {\rm Id}) \ea \ee

\begin{lem} \label{lem:gstar:bif}  We have $g_Q \in L^\infty([0, T], H^{2 s - 2 - d/2}(\R^d)).$ The operator $R_Q$ maps $H^{s'}$ to itself, for any $s' \in \R,$ uniformly in $t \in [0,T].$
\end{lem}

\begin{proof} By regularity of $Q$ and $\d_t Q$ (see the first paragraph of this Section, and Section \ref{sec:spec}), we may follow the proof of Lemma \ref{lem:gstar}.
\end{proof}

\subsection{Weighted Lebesgue spaces in time} \label{sec:weight}

As described in Section \ref{sec:spec}, the real parts of the bifurcating eigenvalues are $O(t).$ Thus the regularization in forward time in the positive elliptic equations will involve an $O(t)$ prefactor. As a consequence, we have to consider weighted Lebesgue spaces in time, and, given $q \in (0,1],$ $T > 0$ and $s' \in \R,$ we let 
\be \label{def:t:space}
 \| f \|_{L^2_{t^q}([0,T], H^{s'})} := \left(\int_0^T t^q \| f(t)\|_{H^{s'}(\R^d)}^2 \, dt \right)^{1/2},
 \ee
 and denote $L^2_{t^q}([0,T], H^{s'}(\R^d))$ the associated vector space. 
 
\begin{lem} \label{lem:interpol} We have, for any $\a > 0,$ any $p \in [0,1),$ and any $\b < \a(p+1)/2,$  the continuous embedding
$$ L^\infty([0, T], H^s) \cap L^2_t([0,T], H^{s + \a}) \subset L^2_{t^p}([0,T], H^{s + \b}).$$
\end{lem}

\begin{proof} Let $u \in  L^\infty([0, T], H^s) \cap L^2_t([0,T], H^{s + \a}).$ We have, with $\theta > 0$ to be specified later: 
$$ \ba \int_0^T \int_{\R^d} t^{p} \langle \xi \rangle^{2(s + \beta)} |\hat u(t,\xi)|^2 \, dt \, d\xi  & \leq \int_{t^{1-p} \langle \xi \rangle^\theta \leq 1} t^{p} \langle \xi \rangle^{2(s + \beta)} |\hat u(t,\xi)|^2 \, dt \, d\xi\\ & + \quad \int_{1 \leq t^{1-p} \langle \xi \rangle^\theta} t^{p} \langle \xi \rangle^{2(s + \beta)} |\hat u(t,\xi)|^2 \, dt \, d\xi.\ea $$ 
We have 
$$ \int_{1 \leq t^{1-p} \langle \xi \rangle^\theta} t^{p} \langle \xi \rangle^{2(s + \beta)} |\hat u(t,\xi)|^2 \, dt \, d\xi \leq \| u \|_{L^2_{t}([0,T], H^{s + \a})}^2,$$
if $\beta + \theta/2 \leq \a.$ 
 Besides, for $\g > 0$ to be specified later, 
$$ \ba  \int_{t^{1-p} \langle \xi \rangle^\theta \leq 1} t^{p} \langle \xi \rangle^{2(s + \beta)} |\hat u(t,\xi)|^2 \, dt \, d\xi & \leq \int \big (t^{1- p} \langle \xi \rangle^\theta \big)^{-\gamma} t^p \langle \xi \rangle^{2(s + \beta)} |\hat u(t,\xi)|^2 \, dt \, d\xi \\ & \leq \left(\int_0^T t^{p - \g(1 - p)} \, dt \right) \| u \|_{L^\infty([0,T], H^s)},\ea $$ 
if $\beta - \g \theta/2  \leq 0.$ The above upper bound is finite if $p - \g(1 - p) > -1.$ 
Optimizing in $\theta$ and $\g$ we end up with the contraint  
$ \b < \frac{\a (p+1)}{2}.$
\end{proof}

\subsection{The positive elliptic equations} \label{sec:+:bif}

System \eqref{ivp:w} is a family of equations that are decoupled to first order. Consider one such equation, for the $k$th coordinate $v_+$ of $v_Q,$ with $k$ chosen so that this equation  is associated with a branch of eigenvalues $\mu$ of $i A$ such that $\Re e \, \mu > 0$ for $t > 0.$  We denote the equation 
\be \label{eq:+:bif} 
 \d_t v_{+} + \op(\psi^\sharp \mu) v_{+} = g_{+} + (R_Q v)_+.
\ee

Equation \eqref{eq:+:bif} is thus one scalar line (one coordinate) from system \eqref{ivp:w}. In particular, the regularity for $g_Q$ and $R_Q v$ stated in Lemma \ref{lem:gstar:bif} extends to $g_+$ and $(R_Q v)_+.$ 

We have for the posited solution $v_+$ of \eqref{eq:+:bif} a forward-in-time regularization result that is analogous to Proposition \ref{prop:+}. The main difference with Proposition \ref{prop:+} is that $\Re e\, \mu_{|t = 0} = 0,$ so that the equation \eqref{eq:+:bif} is degenerate, which will force us to invoke Lemma \ref{lem:interpol}: 

\begin{prop} \label{prop:+:bif}
If the supports of $\psi$ and $\psi^\sharp$ are small enough, and $T > 0$ is small enough, the solution $v_+(t)$ to \eqref{eq:+:bif} belongs to $C^0([0, T], H^{s + \nu_-}(\R^d)) \cap L^2([0,T], H^{s + \nu_+}(\R^d)),$ for any $\nu_- < 1/2$ and any $\nu_+ < 1.$  
\end{prop}

\begin{proof} We follow the steps of the proof of Proposition \ref{prop:+}. 

\medskip

{\it 1. Spectral bound.} By the description of Section \ref{sec:spec}, and choice of the eigenvalue $\mu,$ we have 
\be \label{low:+:bif}
 \Re e \, \mu \geq  t \g |\xi| >0, 
\ee
 in the support of $\psi^\sharp,$ for some $\g > 0,$ for $t \in [0,T_+],$ for some $T_+ > 0.$ 

\medskip

{\it 2. The $L^2$ estimate.} Given $\nu$ as in \eqref{def:rho:bif}, we apply $J D_x^{s + \nu}$ to the equation then perform an $L^2$ estimate. Here $J$ is the regularizing operator $J = (1 + \delta |D|^2)^{-(1 + \nu)/2} \in S^{- (1 + \nu)},$ which depends on $\delta > 0.$ We let $w = J D^{s + \nu} v_+ \in H^{1},$ so that $w \in L^\infty([0,T], H^1).$  We have (as in \eqref{L2:est}) 
\be \label{L2:est:bif} \ba \frac{1}{2} \Re e \, \d_t (\| w \|_{L^2}^2) + t \g \Re e \, ( \op(|\xi|\psi^\sharp) w,  w)_{L^2} & + \Re e \, \big( \op(\psi^\sharp (\mu - t \g|\xi|)) w, w \big)_{L^2} = \Re e \, f,\ea\ee
with
$$ \ba f & = \big(J D_x^{s+\nu} (g_+ + (R_Q v)_+), w)_{L^2} + ([J D^{s + \nu}, \op(\psi^\sharp \mu)] v_+, w)_{L^2} 
\ea$$
Introducing 
$$ 
 \ba \tilde w := R_1(\psi^\sharp J |\xi|^{1 + s + \nu}, \psi^\flat Q) v - R_1(|\xi|, J |\xi|^{s + \nu} \psi^\flat Q) v - |D| R_1(J |\xi|^{s + \nu}, \psi^\flat Q) v,
 \ea
 $$ 
 so that, for fixed $\delta > 0,$ we have $\tilde w \in L^\infty([0,T], H^1).$ We see as in the proof of Proposition \ref{prop:+} that  the left hand side of \eqref{L2:est:bif} is 
\be \label{L2:est:1:bif} \frac12 \Re e \, \d_t (\| w \|_{L^2}^2) + t \g \| |D|^{1/2}  w \|_{L^2}^2 + \g \Re e \, (\tilde w, w )_{L^2} + \Re e \, \big( \op(\psi^\sharp(\mu - t \g |\xi|)) w, w \big)_{L^2}.\ee
 
\medskip

{\it 3. G\r{a}rding's inequality.} By G\r{a}rding's inequality, 
$$  \Re e \, \big( \pdo((\mu - t \g |\xi|) \psi^\sharp) w,w\big)_{L^2} + C( \|(u, \d_x u)\|_{W^{3 + d,\infty}}) \| w \|_{L^2}^2 \geq 0,$$
for some $C > 0$ which depends neither on $\e$ nor on $\delta.$ The key here is that $\mu$ is pointwise as regular in $x$ as $(u, \d_x u).$ Since $u \in L^\infty([0, T],H^s),$ with $s - 1 - d/2 > 3 + d,$ the norm $\| (u, \d_x u)\|_{W^{3 + d,\infty}}$ is finite, uniformly in $t \in [0,T].$ 

Besides, the difference between the para-differential operator $\op(\psi^\sharp(\mu - t |\xi|))$ and the pseudo-differential operator $\pdo(\psi^\sharp(\mu - t |\xi|))$ is bounded as in the proof of Proposition \ref{prop:+}, by regularity of the bifurcating eigenvalue $\mu:$ 
$$ \| ( \op(\psi^\sharp(\mu - t |\xi|)) -  \pdo(\psi^\sharp(\mu - t |\xi|)) ) w \|_{L^2} \lesssim \| w \|_{L^2}.$$

\medskip

{\it 4. Remainder bounds.} Here we bound the two terms in $f$ and the term in $\tilde w.$  
For the first term in $f,$ we introduce $D^{-\nu}/D^\nu$ factors in the $L^2$ scalar product, as in the proof of Proposition \ref{prop:+}. We also introduce $t^{-p}/t^{p}$ factors, with $p \in (0,1/2)$ such that
\be \label{def:p}
 \nu < (2 p + 1) /4 < 1/2,
 \ee
 as we may by definition of $\nu$ \eqref{def:rho:bif}. This gives 
 $$ \big| (\big(J D_x^{s+\nu} (g_+ + R_+ v_Q), w)_{L^2} \big| \lesssim t^{-p} \| D^{-\nu} J D_x^{s+\nu} (g_+ + R_+ v_Q) \|_{L^2} t^{p} \| D^\nu w \|_{L^2},$$
 implying, with Lemma \ref{lem:gstar:bif}:
 $$ \left|\int_0^T (\big(J D_x^{s+\nu} (g_+ + R_+ v_Q), w)_{L^2} \, dt \right| \lesssim \frac{1}{1 - 2 p} \| v \|_{L^\infty H^s} \| w \|_{L^2_{t^{2p}} H^\nu}.$$
Arguing as in the proof of Proposition \ref{prop:+}, and issuing time weights as above, we find that the time integral of the other term in $f,$ and the term in $\tilde w,$ are bounded similarly. 

\medskip

{\it 5. Conclusion.} Summing up, we obtain 
\be \label{est:3:bif} \ba \frac12  \| w \|_{L^2}^2 & + \g \| w \|_{L^2_t H^{1/2}}^2  \\ &  \leq \frac{1}{2} \| w(0)\|_{L^2}^2 + C \big(1 + \| v \|_{L^\infty([0,T], H^s)} \| w \|_{L^2_{t^{2p}}([0,T], H^{\nu})} + \| w \|_{L^2([0,T],L^2)}^2 \big), \ea\ee 
with $p \in (0,1/2)$ defined by \eqref{def:p}, and a constant $C > 0$ which is independent of $\delta.$ By Lemma \ref{lem:interpol} with $s = 0,$ $\a = 1/2$ and $\nu < \b = (2 p + 1)/4,$  we have 
$$ \| w \|_{L^2_{t^{2p}}([0,T], H^{\nu})} \lesssim  \| w \|_{L^2_t H^{1/2}},$$
so that the weighted norm in the upper bound of \eqref{est:3} can be absorbed into the left-hand side. Thus for $T$ small enough, we find 
\be \label{est:4:bif} \| w \|_{L^2}^2 + \frac{\g}{2} \| w \|_{L^2_t H^{1/2}}^2 \leq C \big(1 +  \| w(0)\|_{L^2}^2 + \| v \|_{L^\infty([0,T], H^s)}^2\big),
\ee
with a constant $C > 0$ that is independent of $\delta.$ This implies that the sequence $w = w_\delta$ is bounded in $L^2_t([0,T],H^{1/2}).$ Applying Lemma \ref{lem:interpol} again, this time with $s = 0,$ $\a = 1/2$ and $p = 0,$ we find that $(w_\delta)$ is bounded in $L^2([0,T], H^{\beta}),$ for any $\beta < 1/2.$ As in the proof of Proposition \ref{prop:+}, this implies that $v_+$ belongs to $L^2([0,T], H^{s + \nu + \beta}),$ with $\nu < 1/2$ and $\beta < 1/2.$ The equation in $v_+$ implies that $\d_t v_+ \in L^2([0,T], H^{s + \nu + \beta - 1}),$ hence by interpolation (see for instance Theorem 3 in paragraph 5.9 of \cite{Ev}) $v_+ \in C^0([0,T], H^{s + \nu + \beta - 1/2}).$ Since we may choose $\nu$ and $\beta$ such that $\nu + \beta > 1/2,$ this implies the result. 
\end{proof}

\begin{rem} \label{rem:garding} If in the proof of Proposition {\rm \ref{prop:+:bif}}, we did not have to use G\r{a}rding, or could use G\r{a}rding for a regularized symbol, then we could relax the condition bearing on $s$ in Theorem {\rm \ref{th:transition}}. We see two possible ways around this issue:
\begin{itemize}
\item we could rescale the time, by letting $\tilde v_+(t) := v_+(\sqrt t).$ This would desingularize the elliptic operator, making it possible to use an $\e$-dependent regularization as in the proof of Proposition {\rm \ref{prop:+}}, with $\e > 0$ small and fixed. But then another issue arises: the time rescaling changes $\mu$ into $\mu(\sqrt t)/(2 \sqrt t),$ which is fine for the real part, but leaves us with a large imaginary part. This situation seems to mandate the use of Lagrangian coordinates. We would then need an Egorov lemma, possibly in the form given by A. Boulkhemair in {\rm \cite{Bou}}. 
\item Alternatively, we could try a $\xi$-dependent regularization, that is use the smoothing kernel from the proof of Proposition {\rm \ref{prop:+}}, with $\e = \langle \xi \rangle^{-\kappa},$ for an appropriate choice of $\kappa \in (0,1).$ Here the issue is that this would change the symbols from classical $S^1_{1,0}$ symbols into H\"ormander classes of symbols $S^1_{1,\kappa}.$ For such symbols G\r{a}rding's inequality holds (as discussed for instance in Remark 1.1.28 and Chapter 2 of {\rm \cite{LNT}}), but the para-to-pseudo approximation result that we need has not been written out in this framework as far as we know. 
\end{itemize}
\end{rem}

\subsection{Conclusion} \label{sec:conclusion:bif}

 For the negative elliptic equations in system \eqref{ivp:w}, we can prove a backward-in-time regularization result by a straightforward adaptation of Proposition \ref{prop:-}, based on the proof of Proposition \ref{prop:+:bif}. 
 
 For the hyperbolic equations in system \eqref{ivp:w}, the propagation of regularity is verified as in Section \ref{sec:hyp}.
 
 This proves that $v_Q$ is a bit more regular than the posited solution $u.$ We may repeat the argument around any $\xi \in \S^{d-1}$ (note indeed that Assumption \ref{ass:bif} assumes the same type of time-differentiable bifurcation at every frequency on the sphere), and conclude, as in Section \ref{sec:endgame}, that the posited solution $u$ itself is more regular than it was original assumed. We may repeat the argument, and after a few iterations obtain a contradiction with Theorem 2 from \cite{NT1}, which concludes the proof.

\appendix

\section{Notation and classical results on pseudo-differential operators} \label{sec:symb}

 We denote $C^{k,\theta} S^m$ the class of classical symbols $a(x,\xi),$ possibly matrix-valued, with regularity $C^{k,\theta}$ in $x$ (H\"older regularity) and order $m,$ for $k \in \N$ with $k < \infty,$ $\theta \in (0,1)$ and $m \in \R.$ Given $a \in C^{k} S^m = C^{k,0} S^m,$ we denote $$\| a \|_{m,k,k'} := \sup_{\begin{smallmatrix}  |\a| \leq k \\ |\b| \leq k' \end{smallmatrix}} \big\|  \,\langle \xi \rangle^{|\b| - m} |\d_x^\a \d_\xi^\b a(x,\xi)| \, \big\|_{L^\infty(R^{2d})}.$$ 
 We define similarly $\| a \|_{m,k + \theta,k'}$ for $a \in C^{k,\theta} S^m,$ by taking a H\"older norm for the order-$k$ spatial derivatives of the symbol. 

 The class $C^\infty S^m,$ for $m \in \R,$ is defined as the set of symbols which are smooth and their arguments $(u,\xi) \in \R^N \times \S^{d-1},$ and such that $\langle \xi \rangle^{|\b| - m} |\d_u^\a \d_\xi^\b \s(u,\xi)| \leq C_{\a\b}(|u|),$ for all $(u,\xi) \in \R^N \times \R^d,$ for all $(\a,\b) \in \N^{N \times d},$ for some nondecreasing function $C_{\a\b}.$ 
 
 It is easy to verify (see Lemma A.1 in \cite{NT1}) that if $b \in C^\infty S^m,$ if $u \in H^s$ with $s > d/2,$ then $a(x,\xi) := b(u(x),\xi)$ belongs to $C^{[s - d/2], s - d/2 - [s - d/2]} S^m.$ 
 
 We denote $\pdo(a)$ the pseudo-differential operator associated with symbol $a$ in classical quantization: 
$$ \pdo(a) u := \int_{\R^d} e^{i x \cdot \xi} a(x,\xi) \hat u(\xi) \, d\xi,$$
for $u \in {\mathcal S}(\R^d),$ 
where $\hat u = {\mathcal F}(u)$ is the Fourier transform of $u.$ We denote $\op(a)$ the para-differential operator associated with $a$ in classical quantization, given some admissible cut-off $\phi_{adm}:$ 
$$
 \op(a) := \pdo\Big(\Big(\big( {\mathcal F}^{-1} \phi_{adm}(\cdot,\xi) \big) \star a(\cdot,\xi) \Big)(x)\Big)
$$

 Given $a \in C^0 S^m,$ the para-differential operator $\op_j(a)$ maps $L^2$ to $H^{-m},$ and
 \be \label{action}  \| \op(a) \|_{H^s \to H^{s-m}} \lesssim \| a \|_{m,0, 1 + [d/2]},\ee 
 the implicit constant depending only on dimensions (see for instance Theorem 5.1.15 in \cite{Metivierbis}).

 Given $m_1, m_2 \in \R,$ given $r  > 0$, given $a_1 \in C^r S^{m_1},$ $a_2 \in C^r S^{m_2}$ if $r \in \N,$ or given $a_1 \in C^{[r], r - [r]} S^{m_1},$ $a_2 \in C^{[r], r - [r]} S^{m_2}$ if $r \notin \N,$ we have 
\be \label{composition:para}
 \op(a_1) \op(a_2) = \sum_{0 \leq k  < r } \frac{(-i)^{|\a|}}{\a!} \op\left(  \d_\xi^\a a_1 \d_x^\a a_2\right) + R_r(a_1,a_2),
\ee 
with a remainder $R_r(a_1,a_2)$
which maps $H^{s + m_1 + m_2 - r}$ to $H^s$ with norm 
 \be \label{remainder} \big\| R_r(a_1,a_2) v \big\|_{H^s} \lesssim \Big( \| a_1\|_{m_1,0,m(d,r)} \| a_2 \|_{m_2,r,d} + \| a_1 \|_{m_1,r,m(d,r)} \| a_2 \|_{m_2,0,d}\Big) \| v \|_{H^{s + m_1 + m_2 - r}},\ee
 for some $m(d,r)$ depending on $d$ and $r$ (see for instance Theorem 6.1.4 of \cite{Metivierbis}).

Given $a \in C^1 S^m,$ denoting $\op(a)^\star$ the adjoint of $\op(a)$ with respect to $L^2,$ and $a^\star$ the transpose of the conjugate of matrix $a,$ we have 
 \be \label{adjoint} \op(a)^\star = \op(a^\star) + R_1^{\star}(a),\ee
 with a remainder $R_1^{\star}(a)$ which maps $H^{s + m - 1}$ to $H^s$ with norm 
 \be\label{remainder:adjoint}
  \| R_1^{\star}(a) \|_{H^{s + m - 1} \to H^s} \lesssim \| a \|_{m,1,d_\star}.
 \ee

For all $Q \in C^{3 + d} S^1,$ 
if $\Re e \, Q \geq 0,$ then G\r{a}rding's inequality states that (see for instance in Theorem 4.32 in \cite{Zworski}, and the proof of Theorem 1.1.26 in \cite{Lerner}): 
\be \label{garding0}
 \Re e \, (\pdo(Q) u, u)_{L^2} + \| Q \|_{0, 3 + d, 3 + d} \| u \|_{L^2}^2 \geq 0.
\ee

\end{document}